\documentclass{article}

\usepackage{PRIMEarxiv}

\usepackage[utf8]{inputenc} 
\usepackage[T1]{fontenc}    
\usepackage{hyperref}       
\usepackage{url}            
\usepackage{booktabs}       
\usepackage{amsfonts}       
\usepackage{nicefrac}       
\usepackage{microtype}      
\usepackage{lipsum}
\usepackage{fancyhdr}       
\usepackage{mathrsfs}
\usepackage{amsmath} 
\usepackage{amssymb} 
\usepackage{graphicx}       
\graphicspath{{media/}}     

\title{Comparison of Deep Learning and Deterministic Algorithms for Control Modeling\thanks{Dated: Aug 25, 2022}
}

\author{
  Hanfeng Zhai \\
  Sibley School of Mechanical and Aerospace Engineering, \\
  Cornell University,\\
  Ithaca, NY 14850, USA\\
\And Timothy Sands\thanks{Correspondence: \texttt{tas297@cornell.edu}}\\Department of Mechanical Engineering,\\ Stanford University,\\ Stanford, CA 94305, USA\\ \vspace{5pt}
}

\begin{document}
\maketitle

\begin{abstract}
Controlling nonlinear dynamics arises in various engineering fields. 
We present efforts to model the forced van der Pol system control using physics-informed neural networks (PINN) compared to benchmark methods, including idealized nonlinear feedforward (FF) control, linearized feedback control (FB), and feedforward-plus-feedback combined (C). The aim is to implement circular trajectories in the state space of the van der Pol system. 
A designed benchmark problem is used for testing the behavioral differences of the disparate controllers and then investigating controlled schemes and systems of various extents of nonlinearities. All methods exhibit a short initialization accompanying arbitrary initialization points. The feedforward control successfully converges to the desired trajectory, and PINN executes good controls with higher stochasticity observed for higher-order terms based on the phase portraits. In contrast, linearized feedback control and combined feed-forward plus feedback failed. Varying trajectory amplitudes revealed that feed-forward, linearized feedback control, and combined feed-forward plus feedback control all fail for unity nonlinear damping gain. Traditional control methods display a robust fluctuation for higher-order terms. For some various nonlinearities, PINN failed to implement the desired trajectory instead of becoming "trapped" in the phase of small radius, yet idealized nonlinear feedforward successfully implemented controls. PINN generally exhibits lower relative errors for varying targeted trajectories. However, PINN also shows evidently higher computational burden compared with traditional control theory methods, with at least more than 30 times longer control time compared with benchmark idealized nonlinear feed-forward control. This manuscript proposes a comprehensive comparative study for future controller employment considering deterministic and machine learning approaches.
\end{abstract}



\section{Introduction\label{intro}}

As early as (at least) the late 19th century, scientists made efforts to design and implement control systems to deal with instability, oscillation and various nonlinear and chaotic phenomena \cite{ppl2}.  Maxwell studied valve flow governors \cite{maxwell}, while more recently Cartwright used the van der Pol equation in seismology to model the two plates in a geological fault \cite{Cartwright99}. Fitzhugh \cite{FitzHugh61} and \cite{Nagumo62} used the equation to model action potentials of neurons. Systems exhibiting strong nonlinear behavior are tough problems to control. The standard practice of base controls on the linearization of the system is often rendered ineffective due to the elimination of the nonlinear features. Machine learning is one approach with seeming applicability due to its ability to learn and control nonlinear features.

\subsection{Physics-informed machine learning}

There has been significant recent progress in the field of machine learning in recent decades, starting from the late 80s following the utter failure to achieve its "grandiose objectives" in the 1970s. \cite{Lighthill72} Taking advantage of "big data" and advanced computing technologies such as GPU and TPU computing, there has been exponential growth in the field of deep learning. Central Processing Units (CPU) manage all the functions of a computer and can be augmented by Graphical Processing Units (GPU) and Tensor Processing Units (TPU) to accelerate calculations with application-specific integrated circuits. In the past five years, an explosion of research has re-instantiated "grandiose objectives" manifest in "deep learning". There have been attempts to insert physical information into neural networks (NN) since at least the 1990s, \cite{kbnn1}  relying both on statistical and symbolic learning, called hybrid learning \cite{kbnn2, kbnn3, kbnn4, kbnn5}. Towell, et al. \cite{kbnn2} described hybrid learning methods using theoretical knowledge of a domain and a set of classified examples to develop a method for accurately classifying examples not seen during training. Towell, et al. \cite{kbnn3} introduced methods to refine approximately correct knowledge to be used to determine the structure of an artificial neural network and the weights on its links, thereby making the knowledge accessible for modification by neural learning. Towell, et al. \cite{kbnn4} illustrated a method to efficiently extract symbolic rules from trained neural networks.

Meanwhile, the recent development of physics-informed neural networks (PINNs), originally introduced in 2017 \cite{pinn}, encode differential equations in the losses of the NNs as a soft constraint enabled by automatic differentiation \cite{george}, allowing fast, efficient learning of physical mapping with relatively less labeled data. One well-known application is in the field of fluid fields \cite{science, bubblenet}. An aspect not well known or studied is the implementation of control signals for nonlinear systems using PINNs enabled by inserting the control signals and positional constraints into the loss. This aspect is known as physics-informed deep operator control (PIDOC) \cite{pidoc}. Particularly, it is shown in this work, PIDOC can successfully implement controls to nonlinear van der Pol systems, yet fails to converge to the desired trajectory when the system's nonlinearity is large. 

\subsection{Deterministic algorithms}

In 2017, Cooper et al. \cite{cooper} illustrated how an idealized nonlinear feedforward very effectively controlled highly nonlinear van der Pol systems with fixed parameters, while \cite{pidoc} adopted Cooper's method as the benchmark for comparison, as done here in this manuscript. Based on the work presented in this manuscript on NN-based control and deterministic algorithms, it can be deduced that challenging problems remain open, particularly regarding controlling highly nonlinear systems. The ”grandiose objectives” referred by Sir Lighthill \cite{Lighthill72} remain unfulfilled, and this insight guides both industry and academia efforts in controller design and system stability analysis.

There have also been attempts at comparing classical PID controllers with neural networks \cite{pid_nn_ieee}, refining PID controllers with neural networks \cite{pid_nn_mdpi, pid_nn_mit}. or inserting neural networks into traditional controllers in general \cite{nn_control1, nn_control2, nn_control3}. Hagan et al. \cite{nn_control1} provides a quick overview of neural networks and explains how they can be used in control systems. Nguyen et al. \cite{nn_control2} demonstrated a neural network can learn of its own accord to control a nonlinear dynamic system, while Antsaklis \cite{nn_control3} evaluated whether neural networks can be used to provide better control solutions to old problems or perhaps solutions to control problems that have proved to confound.

Inserting nonlinear approximation by neural networks to refine control and stability is not a new thing and is considered a type of "learning control" dating back to the 80s and 90s. Notwithstanding, as already introduced in \cite{pidoc}, building control frameworks solely with neural networks is relatively rare. Acknowledging the deficiency of related works, this manuscript provides a fairly comprehensive analysis of PIDOC \cite{pidoc} as well as the original methods proposed by Cooper et al. \cite{cooper} on the van der Pol system as a nonlinear representation of oscillating circuits, amongst other example applications. A benchmark is designed considering both the works and analysis of the systematic behavior. Afterward, desired trajectories were modified from the benchmark problem to check how the control methods differ by testing their first and second-order phase portraits.

In this manuscript's Section \ref{sec_problem} we briefly formulate the problem with a brief introduction to the van der Pol system and control schemes. In Section \ref{sec_method} we introduce the control approaches, including physics-informed deep operator control (Section \ref{sec_pidoc}), containing deep learning (Section \ref{sec_dl}) and physics-informed control (Section \ref{sec_pi}); and control theory algorithms in Section \ref{sec_dai}, with linearized feedback control (Section \ref{sec_linear}); idealized nonlinear feed-forward control (Section \ref{sec_nonlinear}) and combined control (Section \ref{sec_combine}); we then briefly how we compare the methods in Section \ref{sec_compare}. Next Section \ref{sec_results} includes results comparing the control schemes: Section \ref{sec_bench} shows how the methods differ on the benchmark problem; Section \ref{sec_trajresult} shows how changing desired trajectories variate the controlled schemes; Section \ref{sec_nonresult} shows how variegating systematic nonlinearity changes different control results. 

\section{Problem Formulation\label{sec_problem}}

\begin{figure}
    \centering
    \includegraphics[scale=0.18]{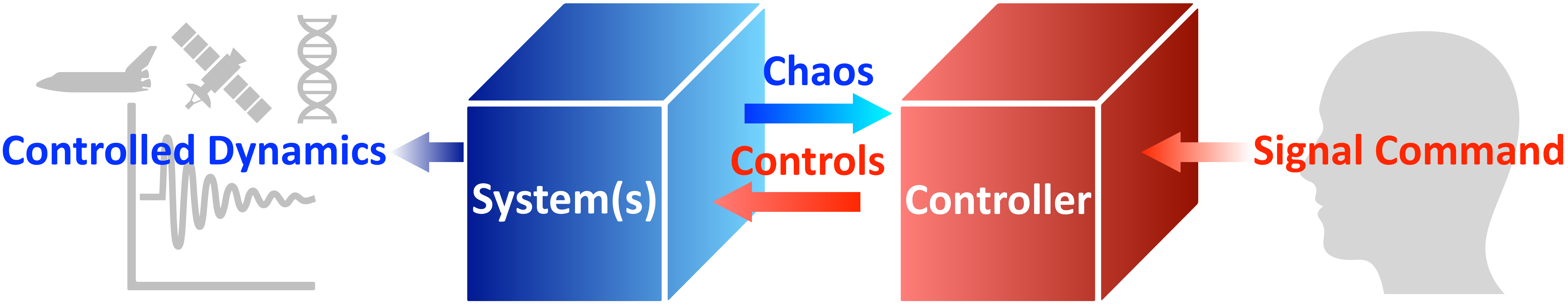}
    \caption{A basic schematic diagram for a control process. The human desired signal command is input to the system through the controller as illustrated in the red box, which passes the control to the targeted system in a "feedforward-feedback-control" loop. Note that the "chaos" from the systems as in the blue box is passed to the controller through the sensor. The final controlled dynamics are output to different applications as illustrated in the left schematic marked as "controlled dynamics". Detailed description please see text.}
    \label{fig_control_schematic}
\end{figure}

As introduced in Section \ref{intro}, the main goal is a comparison of control methods. A basic system schematic is illustrated in Figure \ref{fig_control_schematic}. The command signal was calculated by the controller, passing the control commands to the system, where the system's nonlinear behavior is sensed and fed back using a sensor (not illustrated in the schematic). As the control loop stabilized, the controlled dynamics are output for real-world applications. This manuscript mainly focuses on the controller (red box in figure \ref{fig_control_schematic}), PIDOC, and other control methods are all codified in the controller box.

\begin{figure}
    \centering
    \includegraphics[scale=0.5]{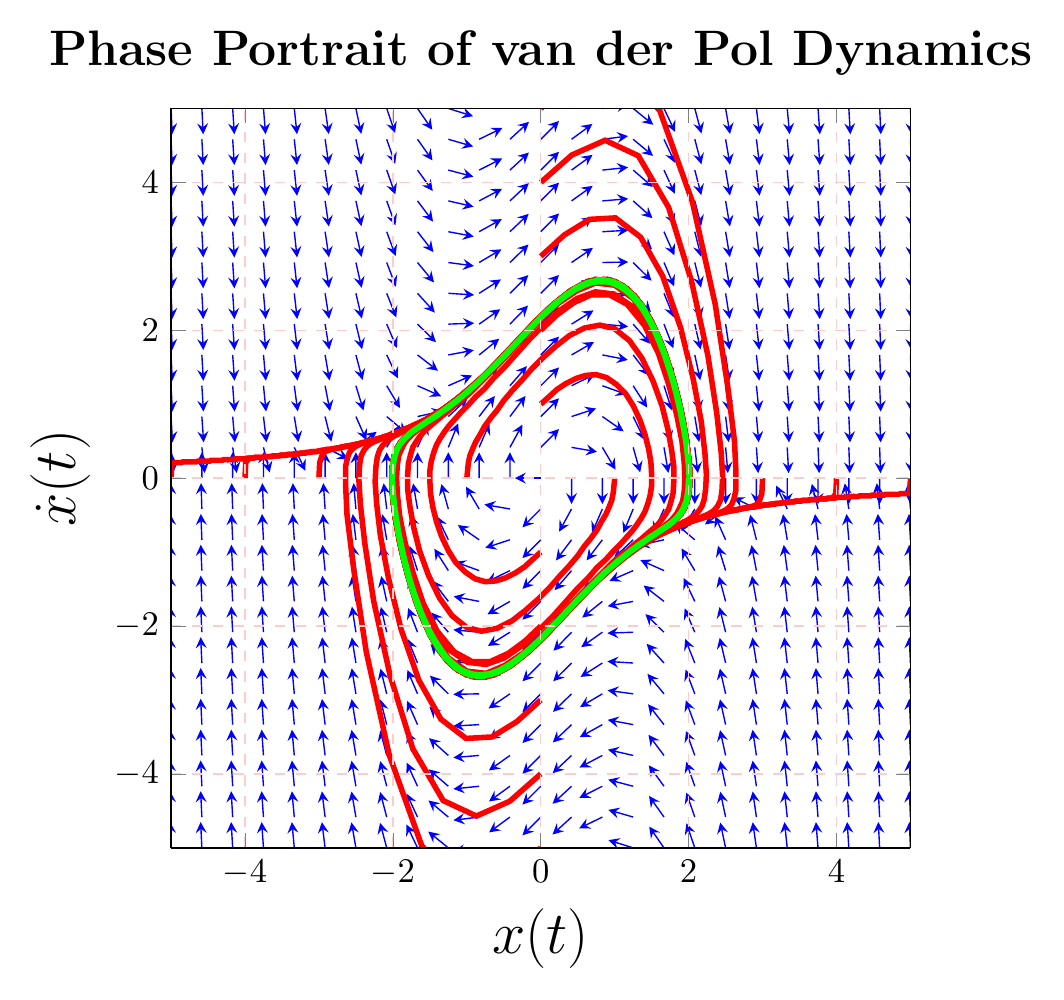}
    \caption{The inherent dynamics of the van der Pol equation \cite{vdp_plot}. The light green line indicates the limit cycle, the manifestation of the strong nonlinearity of the van der Pol inherent dynamics. Disparate red lines indicate trajectories beginning at various initial points, which all eventually fall onto the inherent limit cycle. The blue arrows indicate the total phase field of the van der Pol inherent dynamics, indicating the flow directions.}
    \label{vdp_plot}
\end{figure}

The van der Pol system was adopted to test the control signals' implementation, and a phase portrait of the van der Pol system is illustrated in Figure \ref{vdp_plot} where the system is arbitrarily initialized. Given arbitrary initial points, the trajectory always becomes "entrapped" in a nonlinear track (called a limit cycle), while control methods strive to release the trajectory from the trapped path along the limit cycle and drive the trajectory to some desired, commanded behavior. Such a system was first discovered by van der Pol when investigating oscillating circuits, taking the form \cite{vdp_1,vdp_2}. van der Pol \cite{vdp_1} introduced an oscillatory system with damping that is negative. Together with van der Mark \cite{vdp_2}, he also illustrated how to design an electrical system such that alternating currents or potential differences will occur in the system, having a frequency that is a whole multiple of the forcing function.

\begin{equation}
    \frac{d^2 x}{dt^2} - \mu (1 - x^2) \frac{dx}{dt} + x = 0\label{vdpeq}
\end{equation}
where in the original circuits formulation, $x(t)$ is the current measured in amperes, as the rate of change of the charge \cite{duke} and $\mu$ is a scalar parameter indicating the nonlinearity and the strength of the negative damping \cite{pidoc}. Henceforth, $x(t)$ is referred to as position. 

For testing the proposed methods, control signals are formulated and passed forward to the nonlinear system as commands. The simulated system duplicated the system introduced in \cite{pidoc}, where the MATLAB command \texttt{odeint} solves the equations providing data to feed the training of PIDOC. The van der Pol equation was solved in the time domain from time, $t = [0,50]$, and interpolated with 5000 points. The error control parameters {\tt rtol} and {\tt atol} are $10^{-6}$ and $10^{-10}$, respectively \cite{vdpsolver}.

\section{Methodology and Materials\label{sec_method}}

This section briefly outlines the theoretical foundation of the physics-informed neural network-based algorithm and the alternative based on traditional control theory. The methodology of subsequent numerical experiments used for testing the methods is also introduced. 

\subsection{Physics-Informed Deep Operator Control\label{sec_pidoc}}

\subsubsection{Deep learning\label{sec_dl}}

Physics-informed deep operator control is enabled by the general deep neural network framework, where for the van der Pol system the position is inferred based on the input time domain in accordance with Equation (\ref{nn}).\begin{equation}
    \begin{aligned}
    x_{pred} = (K_L \circ \sigma_L \circ ... \circ K_1 \circ \sigma_1 \circ K_0) t\label{nn}
    \end{aligned}
\end{equation}$K_1, K_2, ..., K_L$, are linear layers; $\sigma_1, \sigma_2, ..., \sigma_L$ are the activation functions, where PIDOC employs \texttt{tanh} activation functions. 

A supervised machine learning framework is defined using external training data, as a formulation minimizing the loss function so that the neural network can capture data features through an optimization process, whereas in traditional neural network approaches $\mathcal{L}$ is usually the difference (errors) between the neural network predictions and training data. Let $\mathcal{L} = \mathcal{L}(t, {\bf p})$ denote the loss function, where $t$ is the input time series and $\bf p$ is the parameter vector contained in formations of $\mathcal{I}$, $\mathcal{D}$, and neural network. Since no external constraints or bounds are enforced, the optimization problem hence takes the form of Equation (\ref{eqn3}) \cite{pidoc}.
\begin{equation}
    \min_{t \subset \mathbb{R}^{d_{out}}} \mathcal{L}(t, {\bf p})
    \label{eqn3}
\end{equation}

Minimizing $\mathcal{L}$ requires reiterating the neural network as defined for the "training". The limited-memory Broyden–Fletcher–Goldfarb–Shanno optimization algorithm, a quasi-Newton method (\textsf{L-BFGS-B} in \texttt{TensorFlow 1.x})  \cite{bfgs, tensorflow} is adopted. Optimization is carried over iterations looping from the blue box (neural network) to purple box ($\mathcal{I}$ \& $\mathcal{D}$) to red box ($\mathcal{L}$) displayed in Figure \ref{schematic}. The maximum iterations are set as $2\times10^5$. In the PIDOC formulation, $\mathcal{L}$ is calculated based on mean square errors of the encoded information to be construed in Section \ref{sec_pi}.

\subsubsection{Physics-informed control\label{sec_pi}}

According to reference, \cite{pidoc}, the control function is enabled by encoding the control signal into the loss function of the neural network, inspired by the formulated physics-informed neural networks (PINNS) \cite{pinn}, where the loss function is computed through the mean squared errors (MSE) elaborated in Equation (\ref{eqn4}). \begin{equation}
    \begin{aligned}
    \mathcal{L} = MSE_{NN} + MSE_{\mathcal{I}} + MSE_{\mathcal{D}}
    \label{eqn4}
    \end{aligned}
\end{equation}where $MSE_{NN}$, $MSE_{\mathcal{I}}$, $MSE_\mathcal{D}$ stands for the neural network generation errors, the initial position loss, and the control signal loss, respectively, computed as Equation (\ref{eqn5}). \begin{equation}
\begin{aligned}
MSE_{NN} &:= \frac{1}{N} \sum_{i=1}^N \left[ x_{train} - x_{pred} \right]^2\\ MSE_{\mathcal{I}} &:= \frac{1}{N} \sum_{i=1}^N \left[ x_{pred}^0 - x_{\mathcal{D}}^0\right]^2\\ MSE_{\mathcal{D}} &:= \frac{1}{N} \sum_{i=1}^N \left[ \left(\frac{d x_{\mathcal{D}}^2}{dt^2} - \frac{d x_{pred}^2}{dt^2}\right) + \left(x_\mathcal{D} - x_{pred}\right)\right]^2
    \label{eqn5}
    \end{aligned}
\end{equation}where $x_{\mathcal{D}}^0$ denotes the initial position of desired trajectory; $x_{pred}$ is the neural network predicted output; $x_{train}$ is the given training data (from system simulation); and $x_{pred}^0$ and $x_{\mathcal{D}}^0$ denote the initial positions of the neural network predicted output and desired trajectory. Detailed formulations are elaborated by reference \cite{pidoc}.

To impose the triangular function signals, we simply impose the form of $x_\mathcal{D}$ in Equation ({\ref{eqn6}}).

\begin{equation}
    x_\mathcal{D}(t) = \Lambda\sin(t),\ \Longrightarrow \dot{x}_\mathcal{D}(t) = \Lambda\cos(t), \ \ddot{x}_\mathcal{D}(t) = -\Lambda\sin(t)
    \label{eqn6}
\end{equation}

Based on such an $x_\mathcal{D}$, the output phase portrait ($\dot{x}(t)$ versus $x(t)$ phase portrait) is expected to be a circular trajectory. To implement different amplitudes of the desired trajectory $\Lambda$, we modify Equation (\ref{eqn5}) to encode the amplitude information into the neural network losses, given the same training data resulting in Equation ({\ref{eqn7}}). 

\begin{equation}
MSE_{NN} := \frac{1}{N} \sum_{i=1}^N \left[ x_{train} - \frac{x_{pred}}{\Lambda}\right]^2
    \label{eqn7}
\end{equation}

where the above equations represent the general formulation of PIDOC. The detailed graphical representation is illustrated as in Figure \ref{schematic} {\bf B}: the control system (deep blue box) first generates nonlinear data that feeds into the neural network, forwards the output to encode the control signals as shown in the deep red box into the loss function through automatic differentiation, and reiterates the training of the neural network until the control signal is fine-tuned for systematic output.

\subsection{Deterministic Control Algorithms\label{sec_dai}}


For the alternative application of control theory, the general framework begins with the modification of Equation (\ref{eqn8}), where controller gains are calculated through the Ricatti equation becoming a controller known as the linear quadratic regulator (LQR) \cite{cooper}.

\begin{equation}
    \frac{d^2 x}{dt^2} - \mu (1 - x^2) \frac{dx}{dt} + x = F (t)
    \label{eqn8}
\end{equation}where $F(t)$ is forced on the nonlinear systems to exert the control. By modifying $F(t)$, different type of controls are implemented, where in our approach we adopt nonlinear feed-forward ($\mathcal{FF}$), linearized feedback control ($\mathcal{FB}$), and the combined controls, to be elaborated in Sections \ref{sec_nonlinear}, \ref{sec_linear}, and \ref{sec_combine}, respectively.

\subsubsection{Linearized feedback control\label{sec_linear}}

In control theory and sciences, a common first step in control design is linearizing nonlinear dynamic equations and then designing the control based on that linearization. For the van der Pol dynamics, Equation (\ref{eqn8}) can be linearized and reduced into Equation (\ref{eqn10}), expressed in state-variable formulation from which state space trajectories are displayed on phase portraits  \cite{cooper}. 




\begin{equation}
    \frac{d\mathscr{X}}{dt} = {\bf A}\mathscr{X} + {\bf B}\mathscr{U}
\label{eqn10}
\end{equation}

The infinite-horizon cost function given by Equation ({\ref{eqn11}})

\begin{equation}
    J = \int_0^{t_{end}} [\mathscr{X}^{\sf T} {\bf Q} \mathscr{X} + \mathscr{U}^{\sf T}{\bf R}\mathscr{U}]dt,\quad {\bf Q} = {\bf Q}^{\sf T}\succeq 0,\ {\bf R} = {\bf R}^{\sf T} \succ 0
\label{eqn11}
\end{equation}

The goal is to find the optimal cost-to-go function $J^* (x)$ which satisfies the Hamilton-Jacobi-Bellman Equation ({\ref{eqn12}})

\begin{equation}
    \begin{aligned}
    \forall \mathscr{X}, \ 0 = \min_u \left[ \mathscr{X}^{\sf T} {\bf Q}\mathscr{X} +  \mathscr{U}^{\sf T} {\bf R}\mathscr{U} + \frac{\partial J^*}{\partial \mathscr{X}} ({\bf A}\mathscr{X} + {\bf B}\mathscr{U})\right]
    \end{aligned}
\label{eqn12}
\end{equation}where to find solutions, Equation ({\ref{eqn13}}) is formed necessitating solution of ({\ref{eqn14}}) which is the {\em algebraic Riccati equation}. The solution of the equation is of well-known form. \footnote{Full derivation: \url{http://underactuated.mit.edu/lqr.html}} Note that the computation of $K_p$, $K_d$, and $[S]$ are based on \textsc{Matlab}\textregistered$\ $command \texttt{[K,S,E] = lqr(A,B,Q,R)}

\begin{equation}
    J^* (\mathscr{X}) = \mathscr{X}^{\sf T}{\bf S}\mathscr{X},\ {\bf S} = {\bf S}^{\sf T}\succeq 0
\label{eqn13}
\end{equation}

\begin{equation}
    0 = \bf SA + A^{\sf T} S - SBR^{\rm -1} B^{\sf T}S + Q
\label{eqn14}
\end{equation}where $\bf A$ and $\bf B$ are the expressions used in the linear-quadratic optimization leading to a feedback controller with linear-quadratic optimal proportional and derivative gains for $K_p$ and $K_d$. The closed loop dynamics are established by Equation (\ref{eqn15}) where the van der Pol forcing function $F(t)$ is a proportional-derivative (PD) controller whose gains used in this manuscript are from \cite{cooper}.


Adopting the linearized feedback control by Cooper et al. \cite{cooper}, Equation (\ref{eqn8}) can thence be expanded in the form:
\begin{equation}
    \begin{aligned}
    \frac{d^2 x}{dt^2} - \mu(1 -x^2) \frac{dx}{dt} + x \equiv F_{\mathcal{FB}}(t) = -K_d (\dot{x}_d - \dot{x}) - K_p (x_d - x)
    \end{aligned}
\label{eqn15}
\end{equation}where $x_d$ is the desired trajectory; $K_d$ and $K_p$ are the derivative and proportional gain, respectively. Similar with our approach in Equation (\ref{eqn6}), $x_d$ is the desired control trajectory, writes $x_d = \Lambda\sin(t)$.




\subsubsection{Nonlinear feed-forward control\label{sec_nonlinear}}
{ 


In idealized nonlinear feed-forward controls, the forced term $F(t) = F_{\mathcal{FF}}(t)$ having the form of the original van der Pol system with the desired trajectory $x = x_d$ executed on:\begin{equation}
    \begin{aligned}
    \frac{d^2 x}{dt^2} - \mu (1 - x^2) \frac{dx}{dt} + x \equiv  F_{\mathcal{FF}} (t) = \frac{d^2 x_d}{dt^2} - \mu (1 - x_d^2) \frac{dx_d}{dt} + x_d
    \end{aligned}
\end{equation} where $x_d$ is the desired signal same as in Equation (\ref{eqn15}). By implementing $x_d$ in the force term, the control is thence applied to the van der Pol system, defined as the nonlinear feed-forward control since the executed force term possesses the form of idealized nonlinear trajectory.

\subsubsection{Combined control\label{sec_combine}}

To apply both the idealized nonlinear feedforward trajectory combined with the linearized feedback, the force term of the combined control simply follows\begin{equation}
    F_{\mathcal{C}}(t) = F_{\mathcal{FF}}(t) + F_{\mathcal{FB}}(t)
\end{equation}}where $F_\mathcal{FB}$ and $F_\mathcal{FF}$ are elaborated in Equations (\ref{eqn14}) and (\ref{eqn15}), respectively. $F_\mathcal{C}$ is then applied to van der Pol system in following the same form as in Equations (\ref{eqn14}) and (\ref{eqn15}).

\begin{figure}[htp]
    \centering
    \includegraphics[scale=0.2]{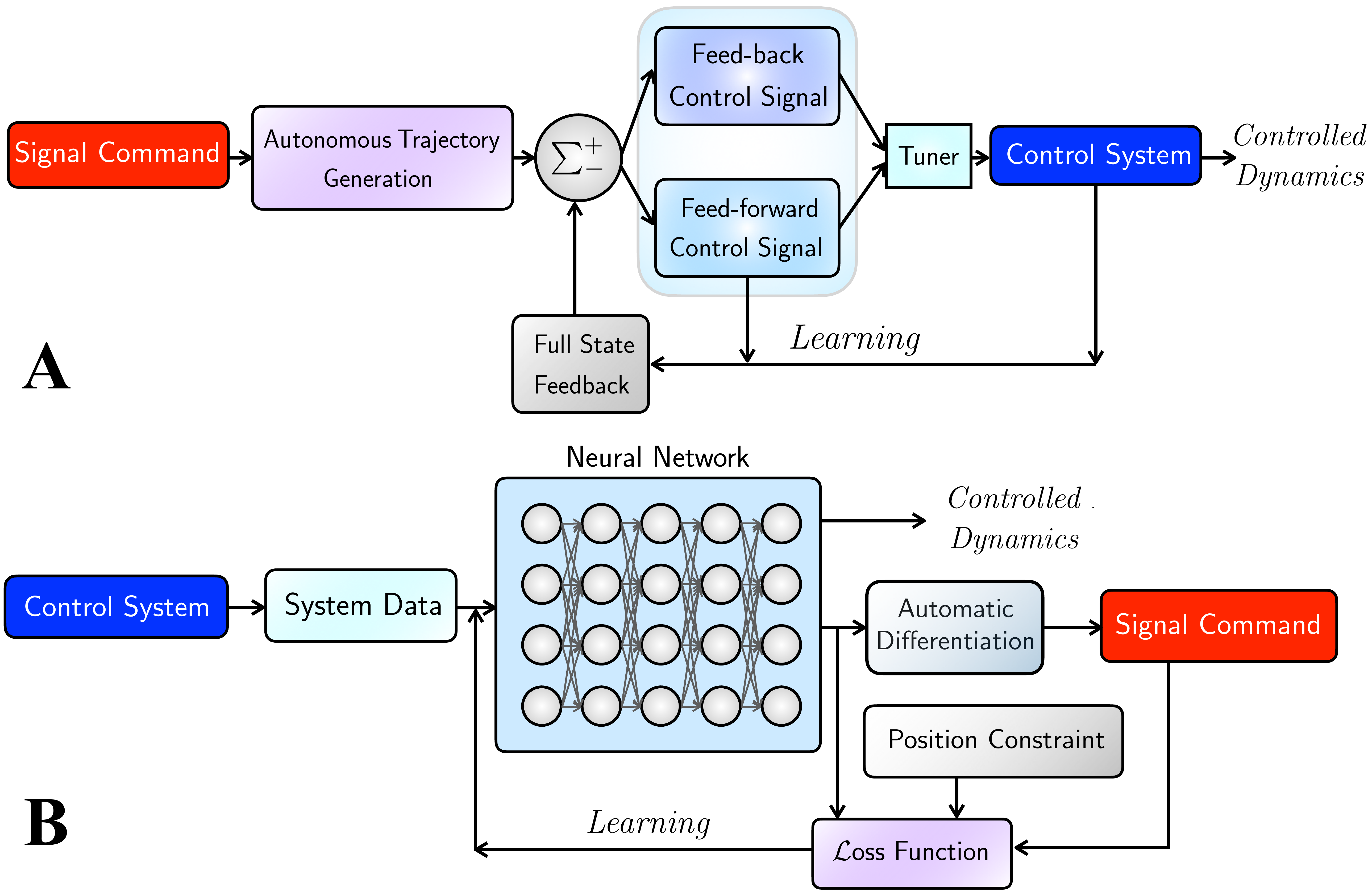}
    \caption{Schematic diagram for the deterministic control algorithms and the deep learning-based PIDOC control scheme. {\bf A}. the schematic for deterministic control algorithms. Note that the light blue tuner can switch the algorithms either to pure idealized nonlinear feed-forward (symbolized as $\mathcal{FF}$, as illustrated in the bottom blue box), linearized feedback (symbolized as $\mathcal{FB}$, as illustrated in the upper dark blue box), or the combined control scheme (symbolized as $\mathcal{C}$, combined both $\mathcal{FF}$ and $\mathcal{FB}$). {\bf B}. the schematic for \textsc{Physics-Informed Deep Operator Control} (PIDOC), symbolized as $\Pi\mathcal{D}$, where the control signal $\mathcal{D}$ (represented in the red box) is inserted in the loss function $\mathcal{L}$ in the purple box as part of the PINN. Detailed description please see the text.}
    \label{schematic}
\end{figure}

{
The basic framework of the controls is shown in Figure \ref{schematic} {\bf A}: the signal command as shown in the deep red box ($x_d$ in our equations) is the first input to the automatic trajectory generator that is forwarded to the gains, and then forwarded to either feed-forward controls ($F_\mathcal{FF}$) on the lower light blue box or feedback controls ($F_\mathcal{FB}$) on the upper dark blue box or the combined approach. The control signals are tuned through the light blue tuner box on the right, which controls the force term applied to the nonlinear system as indicated in the solid blue box on the right. After exerting the desired control signals, the output signals are first fed to the gains as full state feedback indicated in the gray box; the final controlled dynamics are output after the workflow is executed iteratively.
}

\subsection{Comparison and Estimation\label{sec_compare}}


To conduct a fair, decent, and comprehensive comparison of the proposed methods, we consider {\em Systematic analysis} of the provided benchmark problem as we mentioned in Section \ref{intro}, {\em Trajectory convergence} for different amplitudes of desired trajectories, signified by $\Lambda$ in Equation (\ref{eqn6}) and {\em Non-linearity} of the systems with different nonlinearities, signified through $\mu$ ion Equation (\ref{vdpeq}). For the benchmark systematic behavior analysis, considering both the work of Zhai \& Sands \cite{pidoc} and Cooper et al. \cite{cooper}, we pick $\Lambda = 5, \mu = 1$, as a system with low nonlinearity; in which for the PIDOC framework the NN has the structure of $6\times30$. The initial point is picked as (1,0). For systems of different desired amplitudes, $\Lambda$ is changed from $1, 3, 5, 7, 9$. For systems of different non-linearities, $\mu$ is changed from $1,3,5,7,9,10$. The PIDOC was conducted in Google Colab \cite{colab} using \texttt{Python 3.6} compiling \texttt{TensorFlow 1.x} \cite{tensorflow}. Both $\mathcal{FF}$, $\mathcal{FB}$ and $\mathcal{C}$ were written in \textsc{Matlab R2021a} and executed with {\sf Simulink}.

\section{Results and Discussion\label{sec_results}}

\subsection{Benchmark analysis\label{sec_bench}}

\begin{figure}[htp]
    \centering
    \includegraphics[scale=0.18]{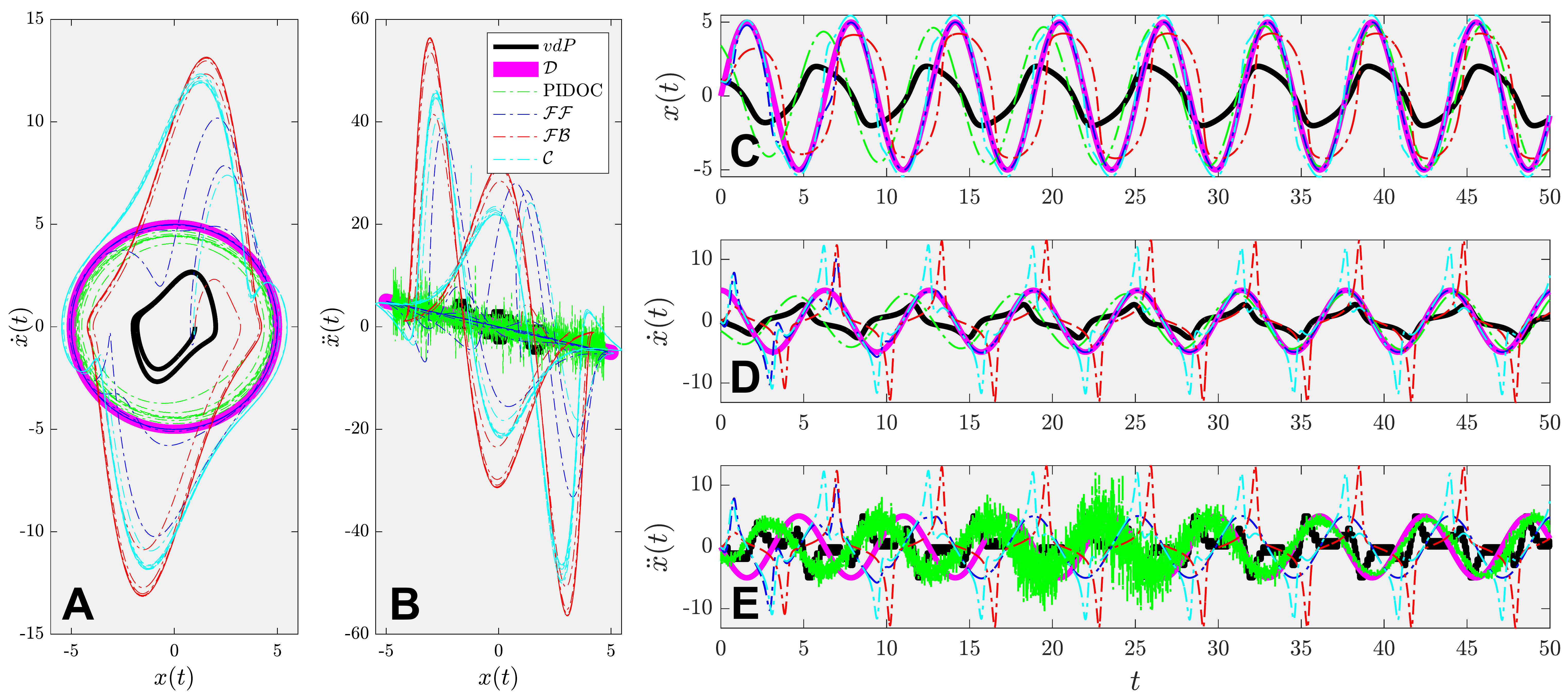}
    \caption{System behavior analysis for the benchmark problem. Note the inherent van der Pol dynamics ($vdP$) is marked in black solid line; the desired trajectory ($\mathcal{D}$) is marked in pink solid line; the PIDOC control is marked in light green dashed line; the feed-forward control ($\mathcal{FF}$) is marked in dark blue dashed line; the feedback control ($\mathcal{FB}$) is marked in red dashed line; the combined control ($\mathcal{C}$) is marked in the light blue line. {\bf A}. the phase portrait of the van der Pol systems of inherent dynamics, desired trajectory, and different control schemes marked in different colors. {\bf B}. the acceleration-position plot. {\bf C}. the time evolution of positions. {\bf D}. the time evolution of velocities. {\bf E}. the time evolution of accelerations.}
    \label{fig_benchmark}
\end{figure}

The results of the benchmark analysis are shown in Figure \ref{fig_benchmark}, where sub-figures {\bf A}, {\bf B} stand for the first and second order phase portraits of different controlled schemes by PIDOC, $\mathcal{FF}$, $\mathcal{FB}$, \& $\mathcal{C}$, marked in different colors dashed lines as elaborated in the caption; compared with the inherent van der Pol dynamics and desired trajectory marked in black and pink solid lines, respectively. Sub-figures {\bf C} to {\bf D} illustrated the time evolution of the zeroth, first, and second order derivatives of the position $x(t)$, with the same color representations as in {\bf A} and {\bf B}. Given the benchmark problem, it can be deduced that all the control theory methods exhibit strong fluctuations at the initial stage of controls, where $\mathcal{FF}$ converge to the desired trajectory successfully as indicated in the deep blue dashed lines, whereas both $\mathcal{FB}$ and $\mathcal{C}$ fails. Another interesting point to be noted is that all the traditional control algorithms exhibit a stronger fluctuation for higher order terms at the beginning stage, yet $\mathcal{FF}$ successfully converge to the trajectory that exhibits better control effects than PIDOC, but $\mathcal{FB}$ and $\mathcal{C}$ displays such a robust fluctuation along the time. To this phenomenon, we provide the following explanation: the errors generated by the linearization of the van der Pol equation accumulate and cause the robust fluctuations as indicated in Figure \ref{fig_benchmark} for the light blue and red lines. However, admittedly, $\mathcal{FF}$ successfully implements the control with higher accuracy for higher order terms than PIDOC; but noted that since $\mathcal{FF}$ only forwarding control signals can be considered as an open-loop system, in real-world practice, trivial noises will be accumulated that leads to the in-feasibility of $\mathcal{FF}$.

\subsection{Trajectory amplitude\label{sec_trajresult}}

\begin{figure}[htp]
    \centering
    \includegraphics[scale=0.3]{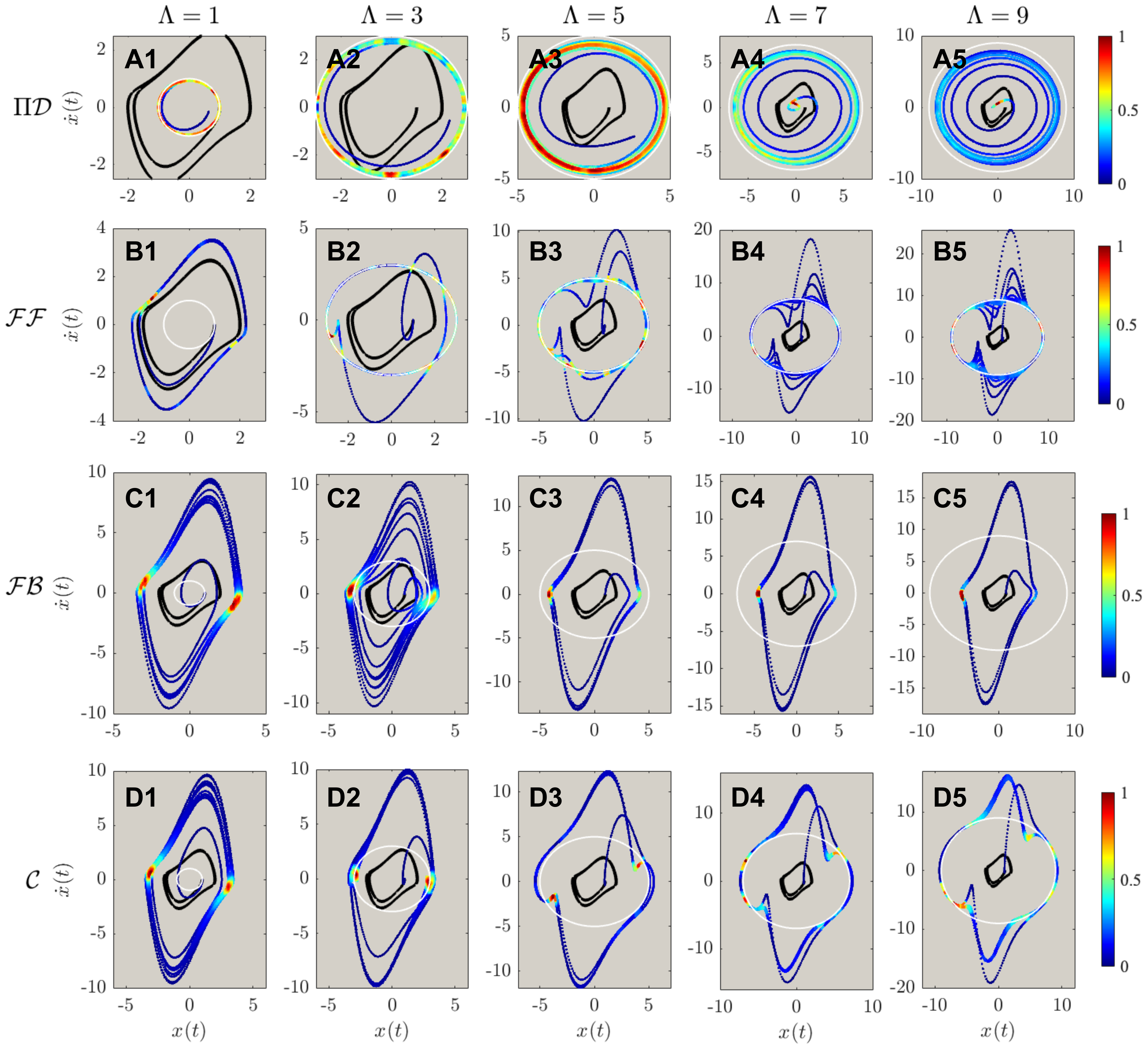}
    \caption{The phase portrait of different controlled schemes, including the inherent van der Pol dynamics, marked in black solid line; the desired trajectory marked in white dashed line; and the controlled dynamics marked in the contoured line. The contour legend was marked from 0 to 1, showing the intensity of the controls. {\bf A1} to {\bf A5} shows the controls by physics-informed deep operator control, symbolized by $\Pi\mathcal{D}$, of different desired trajectories from $\Lambda = 1, 3, ..., 9$. {\bf B1} to {\bf B5} shows the controls by feed-forward controls ($\mathcal{FF}$), from $\Lambda = 1, 3, ..., 9$. {\bf C1} to {\bf C5} shows the controls by feedback controls ($\mathcal{FB}$), from $\Lambda = 1, 3, ..., 9$. {\bf D1} to {\bf D5} shows the controls by feed-forward - feedback combined controls ($\mathcal{C}$), from $\Lambda = 1, 3, ..., 9$. }
    \label{fig1}
\end{figure}

\begin{figure}[htp]
    \centering
    \includegraphics[scale=0.3]{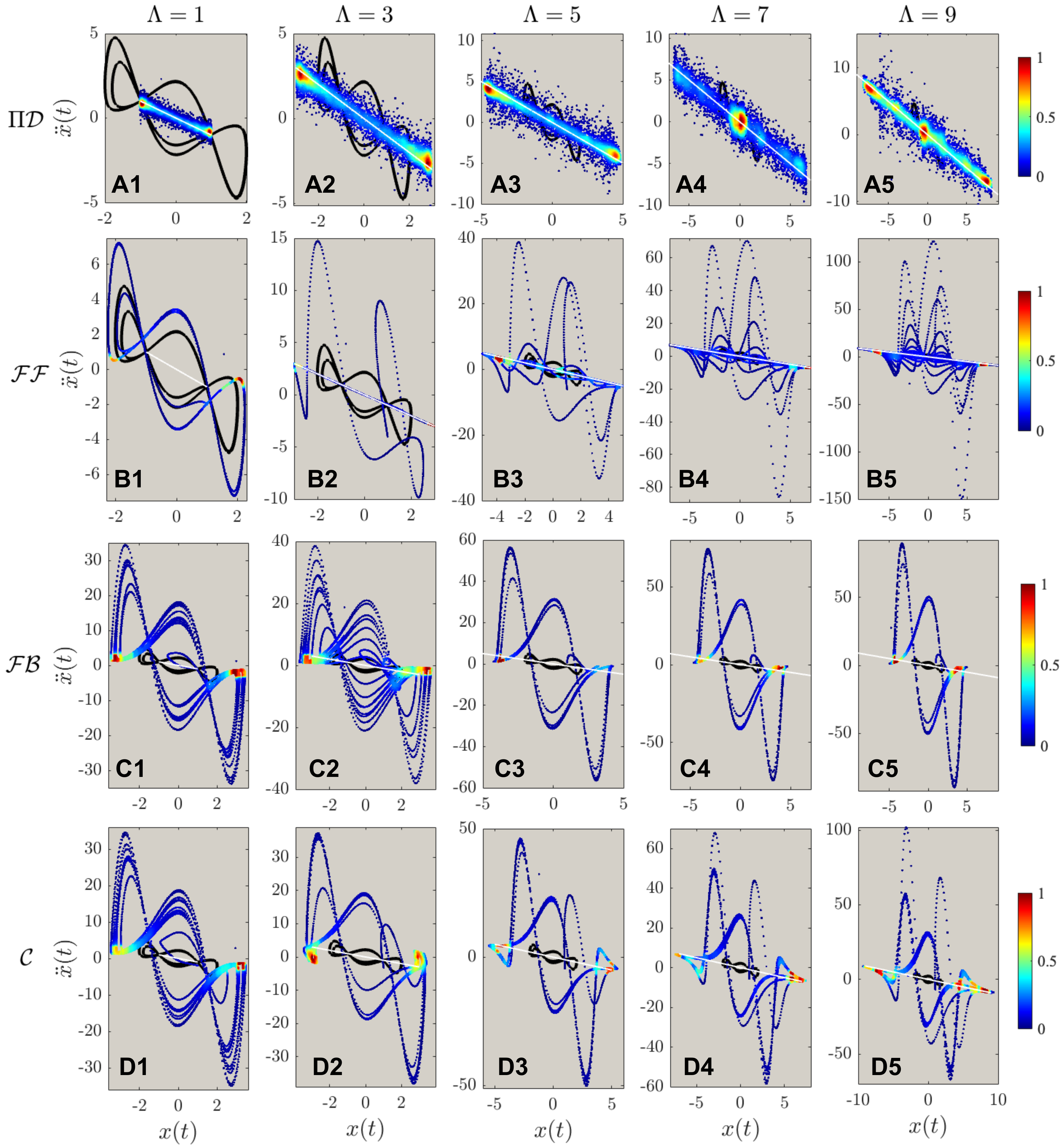}
    \caption{The acceleration-position portrait of different controlled schemes, with the marking colors same as in Figure \ref{fig1}. Note that {$\bf A(1,2,..., 5)$} to {$\bf D(1,2,..., 5)$} are the same as in Figure \ref{fig1}: the implementations of $\Pi\mathcal{D}$, $\mathcal{FF}$, $\mathcal{FB}$, \& $\mathcal{C}$ to different targeted trajectory amplitudes of $\Lambda = 1,3, ..., 9$. }
    \label{fig2}
\end{figure}

The results of controlled dynamics of trajectories of the first and second-order phase portraits are shown in Figures \ref{fig1} and \ref{fig2}, respectively. It can be discerned from Figure \ref{fig1} {\bf A} and {\bf B} that both PIDOC (symbolized as $\Pi\mathcal{D}$ in the figure) and $\mathcal{FF}$ are able to implement controls with an exception of {\bf B1} that $\mathcal{FF}$ failed to control the system when $\Lambda = 1$. Similar with the benchmark problem that both $\mathcal{FB}$ and $\mathcal{C}$ failed to implement the controls with a highly fluctuating behavior, in Figure \ref{fig1} {\bf C} and {\bf D}. An interesting phenomenon reported from {\bf D1} to {\bf D5} is that with increasing trajectory amplitudes we report a better convergence for the combined ($\mathcal{C}$) control. We can hence propose the discussion on such phenomena that for higher values of desired trajectory amplitudes the linearization effect of the feedback reduces for the van der Pol systems.

Figure \ref{fig2} reports the second order phase portraits (acceleration - position diagram) comparing the four methods. Figure \ref{fig2} {\bf A} reports the stochastic approximation nature of PIDOC: the learning-based control executes control signals based randomized sampling for trajectory convergence. Corresponds to Figure \ref{fig1}, {\bf B1} shows the failure of $\mathcal{FF}$ control when $\Lambda = 1$; whereas {\bf B2} to {\bf B5} shows how the second order phase portraits display a higher fluctuation, as also shown from in Figure \ref{fig2} {\bf C} and {\bf D}. Figure \ref{fig2} {\bf B} also shows a strong discretized form of $\mathcal{FF}$ control, as illustrated based on the sparse points. The control contour from both Figures \ref{fig1} and \ref{fig2} both $\mathcal{FB}$ and $\mathcal{C}$ controls (sub-figure \textbf{C} and \textbf{D}) shows an increased control density on the horizontal edges ($x(t)$ direction), indicated by the denser points.

\begin{table}[htp]
    \centering
    \begin{tabular}{c|c c c c c}
        $\Lambda$  & 1 & 3 & 5 & 7 & 9 \\ine
        $\Pi\mathcal{D}$ & 329.82 & 618.97 & 713.30 & 261.37 & 480.77\\
        $\mathcal{FF}$ & 5.70 & 6.01 & 5.26 & 6.21 & 5.66 \\
        $\mathcal{FB}$ & 6.68 & 7.23 & 5.99 & 6.07 & 6.87 \\
        $\mathcal{C}$ & 8.35 & 6.63 & 7.62 & 6.46 & 6.06\\ine
    \end{tabular}
    \caption{The computational burden of the four different control frameworks considering different desired trajectories (desired radius $\Lambda$). Note that the unit are in seconds (default of \texttt{timeit.time()} in \texttt{Python} and \texttt{cputime} in \textsc{Matlab}). }
    \label{tab1}
\end{table}

The total computational burden of the four methods is shown in Table \ref{tab1}: the PIDOC framework shows an evidently larger computing time than $\mathcal{FF}$, $\mathcal{FB}$ and $\mathcal{C}$; generally, $\mathcal{FF}$ execute the fastest control and $\mathcal{C}$ exhibits the longest control time within the tested control theory algorithms. We provide the following explanations for the above phenomena: (1) the PIDOC framework is based on the training of the NN, where the approximation of nonlinear data takes exponentially longer compared with just implementing the control commands; (2) since $\mathcal{FF}$ can be considered as an open-loop implementation of control signals, where the elimination of feedback and error adjustment reducing computation time; (3) the combination of both feed-forward and feedback requires estimation of the route execution and linearizations, consumes more time. Based on the computation time we can discern that although more stable control implementations are exhibited by PIDOC, the drawback is also evident: the considerably longer training time required for implementing the control. 

\subsection{Nonlinear effects\label{sec_nonresult}}

\begin{figure}[htp]
    \centering
    \includegraphics[scale=0.32]{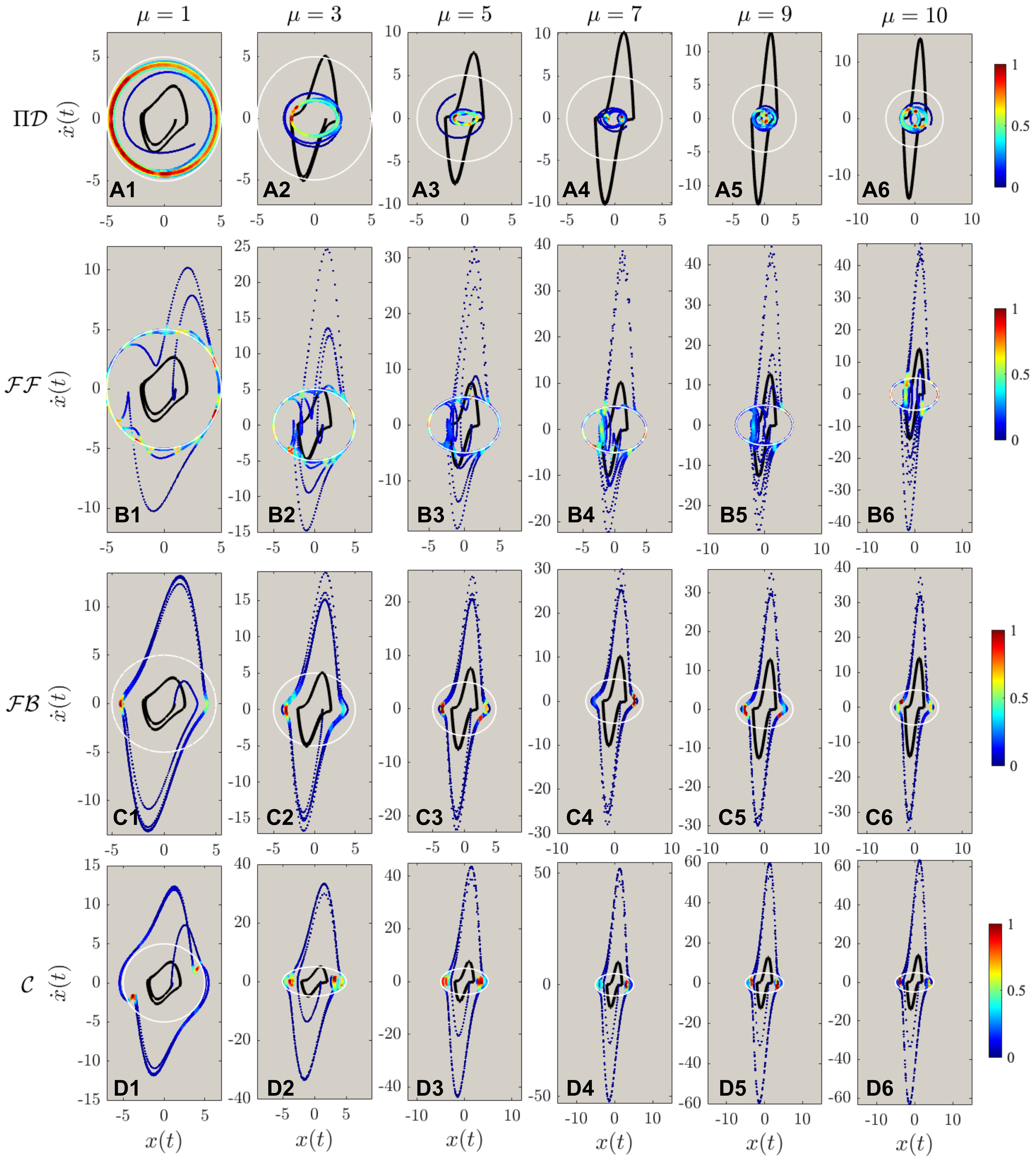}
    \caption{The phase portrait of different controlled schemes, with the marking colors same as in Figure \ref{fig1}. {\bf A1} to {\bf A6} shows the controls by physics-informed deep operator control, symbolized by $\Pi\mathcal{D}$, of different van der Pol systems with different nonlinearities from $\mu = 1,3,5,7,9,10$. {\bf B1} to {\bf B6} shows the controls by feed-forward controls ($\mathcal{FF}$), from $\mu = 1,3,5,7,9,10$. {\bf C1} to {\bf C6} shows the controls by feedback controls ($\mathcal{FB}$), from $\mu = 1,3,5,7,9,10$. {\bf D1} to {\bf D6} shows the controls by feed-forward - feedback combined controls ($\mathcal{C}$), from $\mu = 1,3,5,7,9,10$.}
    \label{fig3}
\end{figure}

The results of different control for systems of different nonlinearities with a fixed desired trajectory $\Lambda = 5$ are shown in Figure \ref{fig3}. Same as reported by Zhai \& Sands \cite{pidoc}, the PIDOC control fails to implement control for systems of high nonlinearities as to be "trapped" in a smaller radius trajectory. The $\mathcal{FF}$ control was implemented successfully, with a strong fluctuation reported for high nonlinearities observed from {\bf B1} to {\bf B5}, with the failed implementation when $\mu = 10$ as shown in {\bf B6}, which can be considered as nonlinearity threshold. Both $\mathcal{FB}$ and $\mathcal{C}$ also failed for control execution same as in Figures \ref{fig1} and \ref{fig2}. To note, both the control theory methods implemented show an evident higher data density along the horizontal edges, which can be adopted to infer the nature of control theory methods: stronger control imposition near edges, corresponding to the wave crests and troughs as for the time evolution of the position.

\begin{figure}[htp]
    \centering
    \includegraphics[scale=0.32]{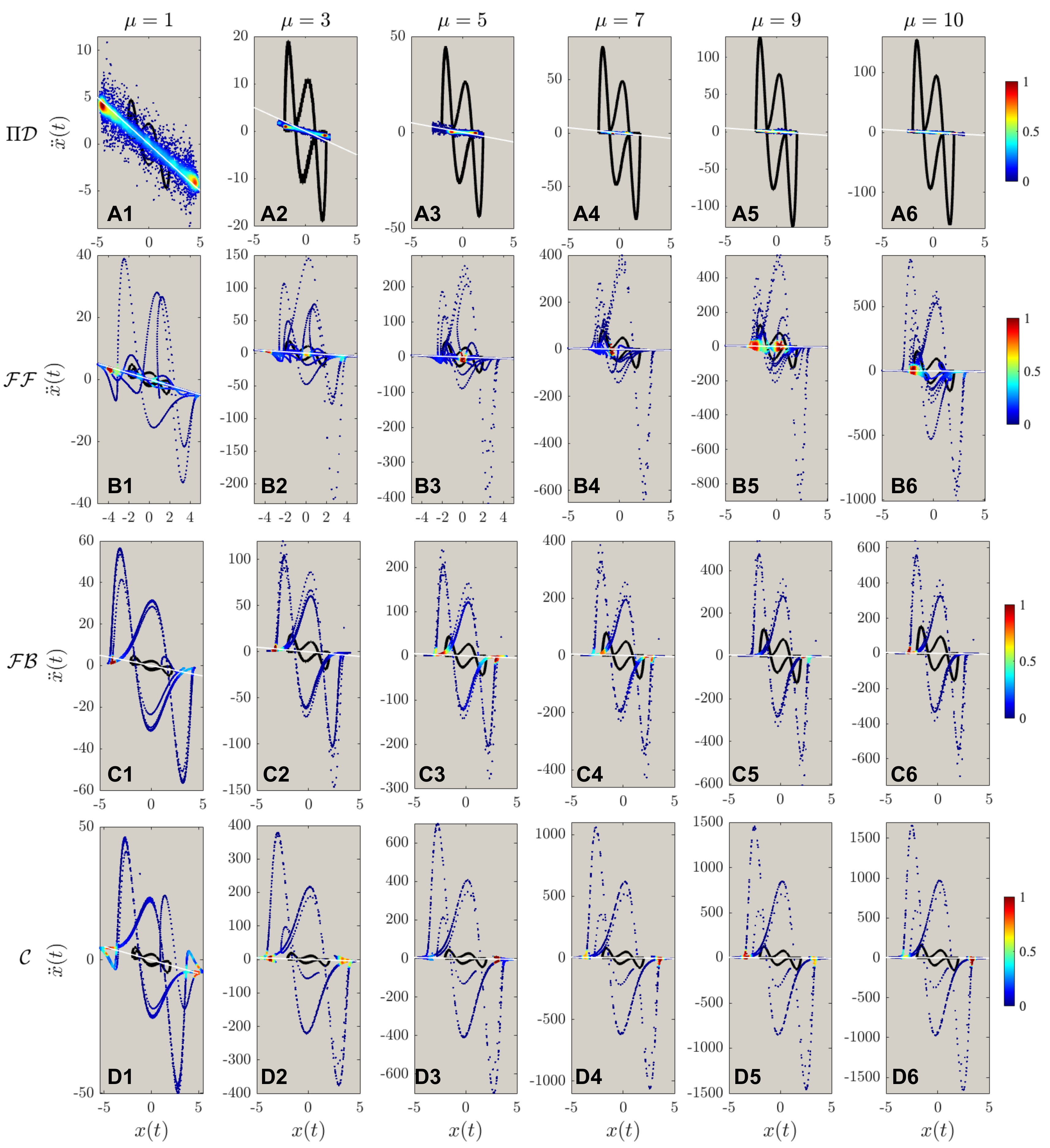}
    \caption{The acceleration-position portrait of different controlled schemes, with the marking colors same as in Figure \ref{fig1}. Note that {$\bf A(1,2,..., 6)$} to {$\bf D(1,2,..., 6)$} are the same as in Figure \ref{fig1}: the implementations of $\Pi\mathcal{D}$, $\mathcal{FF}$, $\mathcal{FB}$, \& $\mathcal{C}$ to different targeted trajectory amplitudes of $\Lambda = 1,3, ..., 9$.}
    \label{fig4}
\end{figure}

\begin{figure}[htp]
    \centering
    \includegraphics[scale=0.35]{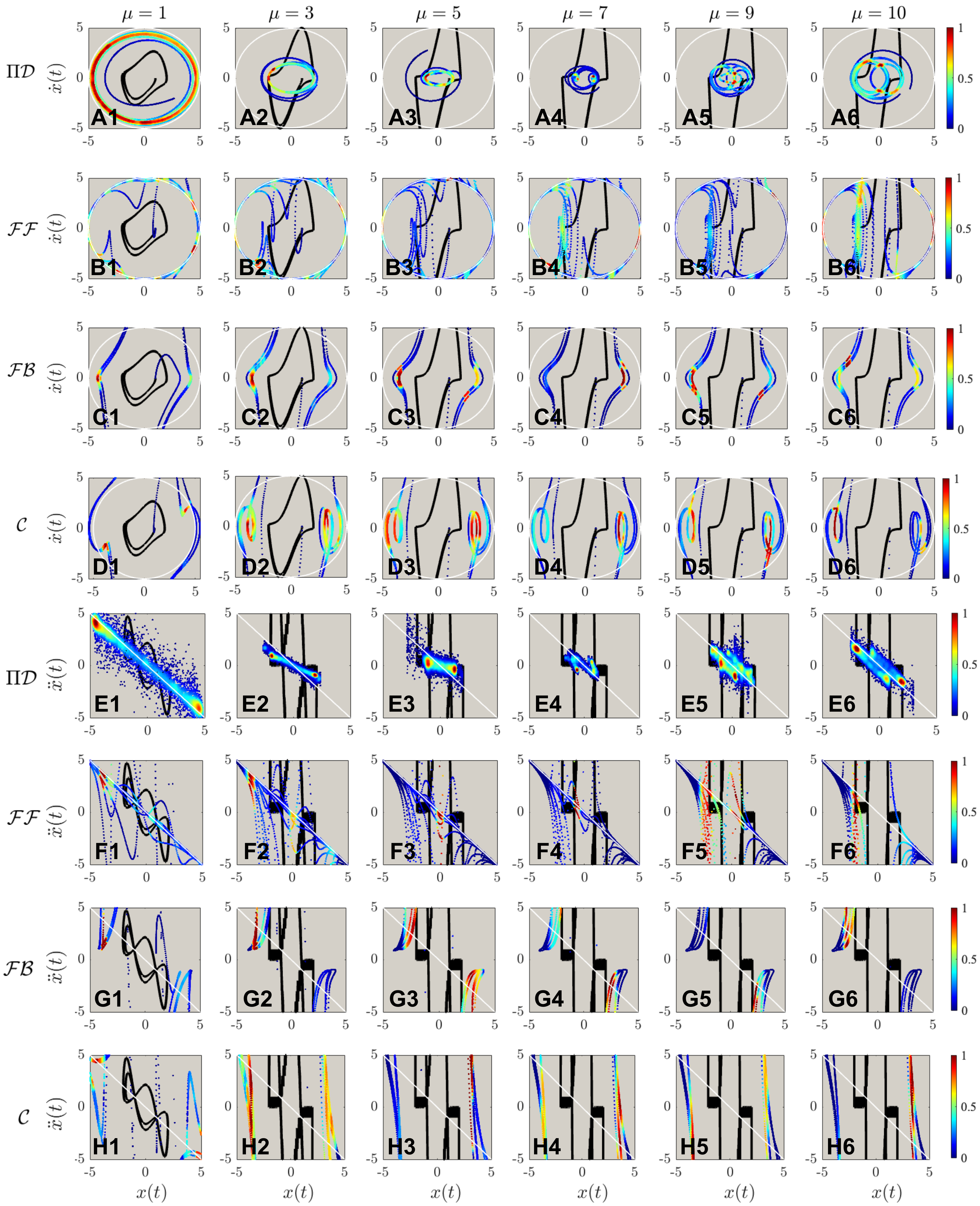}
    \caption{The zoomed view of the controlled schemes portrait corresponding to both Figures \ref{fig3} and \ref{fig4}. Note that the sub-figures {$\bf A(1,2,..., 6)$} to {$\bf D(1,2,..., 6)$} are the same as in Figure \ref{fig3}; whereas the sub-figures {$\bf E(1,2,..., 6)$} to {$\bf H(1,2,..., 6)$} corresponds to sub-figures {$\bf A(1,2,..., 6)$} to {$\bf D(1,2,..., 6)$} in Figure \ref{fig4}. }
    \label{fig5}
\end{figure}

The second order phase portraits are shown in Figure \ref{fig4}: as for the control theory methods, evidently higher nonlinearities are observed for $\mathcal{C}$ compared with $\mathcal{FF}$ and $\mathcal{FB}$; a more discrete points distribution indicate larger steps for control implementations. Just by observing Figure \ref{fig4} {\bf A}, we discern that the systematic nonlinearity was very high as indicated in the black solid line compared with the white dashed line as for the desired trajectory. However, comparing {\bf B} to {\bf D} we observe that for systems of higher nonlinearities, the control displays extremely strong fluctuations at the beginning stage of the control. Based on such a phenomenon we hence deduce another finding for control theory properties: the control implementation will enlarge the nonlinear signals with a more larger steps of control discretization. To present a more detailed analysis of Figures \ref{fig3} and \ref{fig4}, we create Figure \ref{fig5} for a zoomed view of the control schemes for both first and second-order phase portraits. Interestingly, vortex-liked structures are observed in the first-order phase portrait for both PIDOC and $\mathcal{C}$ along the edges of the circular trajectory. Figure \ref{fig5} {\bf B6} shows how $\mathcal{FF}$ fails control imposition in detail: an oscillation along the circular causes the "split" of the controlled trajectory vertically, where such a trend has already been observed in Figure \ref{fig5} {\bf B5}. Figure \ref{fig5} {\bf C} clarifies a phenomenon that has already been observed and discussed: an increased data density along the edges of the desired control schemes indicates a stronger control implementation along the edges.

\begin{table}[htp]
    \centering
    \begin{tabular}{c|c c c c c c}
        $\mu$  & 1 & 3 & 5 & 7 & 9 & 10\\ine
        $\Pi\mathcal{D}$ & 713.30 & 257.64 & 305.91 & 225.15 & 197.76 & 199.52 \\
        $\mathcal{FF}$ & 5.26 & 10.33 & 7.77 & 6.45 & 5.38 & 5.12\\
        $\mathcal{FB}$ & 5.99 & 5.61 & 6.84 & 6.02 & 5.52 & 5.30\\
        $\mathcal{C}$ & 7.62 & 5.94 & 5.12 & 5.83 & 5.26 & 5.42\\ine
    \end{tabular}
    \caption{The computational burden of the four different control frameworks considering different systems of nonlinearities (different $\mu$ values). Note that the unit is in seconds (same as in Table \ref{tab1}). }
    \label{tab2}
\end{table}

The computational burden as shown in Table \ref{tab2} displays similar trends as in Table \ref{tab1}: PIDOC exhibits an evidently higher computation time, attributed to the NN training. $\mathcal{C}$ exhibits a higher control time than $\mathcal{FF}$ and $\mathcal{FB}$. Another interesting phenomenon is: that with the increasing system nonlinearity, PIDOC shows a  decreasing computation time. Corresponds to Figures \ref{fig3}, \ref{fig4}, and \ref{fig5}, we propose the following explanation: as the PIDOC-controlled schemes are entrapped in a trajectory with a lower radius, the NN straining stops at an earlier stage since the optimizer (\textsf{L-BFGS-B}) "discern" that more iterations won't keep decreasing the loss, which leads to a lower computation time but lower quality control. To better quantify the computational burden differences, we create Table \ref{relative_time}, taking nonlinear feed-forward control employed by Cooper et al. \cite{cooper} as a benchmark: PIDOC displays evidently higher computational burden compared with $\mathcal{FF}$, with at least more than 30 times of the benchmark time to up to 100 plus more times.

\begin{table}[htbp]
    \centering
    \begin{tabular}{c | c c c c}
        $\hat{\mathcal{T}}$ & $\Pi\mathcal{\mathcal{D}}$ & $\mathcal{FF}$ & $\mathcal{FB}$ & $\mathcal{C}$ \\ine
        $\Lambda = 1$ & 39.4999  &  1.0000  &  0.8000   & 0.6826\\
        $\Lambda = 3$ & 93.3585  &  1.0000  &  1.0905   & 0.9065 \\
        $\Lambda = 5$ & 93.6090  &  1.0000  &  0.7861   & 0.6903\\
        $\Lambda = 7$ & 40.4602  &  1.0000  &  0.9396   & 0.9613 \\
        $\Lambda = 9$ & 79.3355  &  1.0000  &  1.1337   & 0.9340 \\ine
        $\mu = 1$ & 135.6084  &  1.0000  &  1.1388  &  1.4487 \\
        $\mu = 3$ & 51.5006  &  1.0000  &  0.9444  &  1.7391 \\
        $\mu = 5$ & 38.6258  &  1.0000  &  1.3359   & 1.5176 \\
        $\mu = 7$ & 34.2228   & 1.0000   & 1.0326 &   1.1063 \\
        $\mu = 9$ & 42.8036  &  1.0000  &  1.0494   & 1.0228\\
        $\mu = 10$ & 47.5353  &  1.0000  &  0.9779   & 0.9446\\ine
    \end{tabular}
    \caption{The relative computation time comparing PIDOC and control theory algorithms regarding different trajectories and nonlinearities.}
    \label{relative_time}
\end{table}

To quantify the control errors, we generate Table \ref{relative_error} to compare the control qualities based on the absolute errors. The equation for computing the average absolute relative errors of different control signals are\begin{equation}
    \begin{aligned}
    \left\|\hat{\mathcal{E}}\right\| = \sum_{i = 1}^M \frac{1}{M}\left\|\frac{x_{control} - x_{\mathcal{D}}}{x_{\mathcal{D}}}\right\|
    \end{aligned}\label{relerr}
\end{equation}

It can be observed from Table \ref{relative_error} that PIDOC generally exhibits lower control errors compared with traditional control methods in different trajectories. For different nonlinearities, corresponding to Figure \ref{fig5}, it can be observed that nonlinear idealized feed-forward control exhibits better control qualities.  

\begin{table}[htbp]
    \centering
    \begin{tabular}{c | c c c c}
        $\|\hat{\mathcal{E}}\|$ & $\mathcal{C}$ & $\mathcal{FF}$ & $\mathcal{FB}$ &  $\Pi\mathcal{\mathcal{D}}$\\ine
        $\Lambda = 1$ & 2.1379  &  1.7199  &  2.0618   & 0.2225\\
        $\Lambda = 3$ & 0.3645  &  0.4124  &  0.4473   & 0.2102\\
        $\Lambda = 5$ & 0.3884  &  0.4245  &  0.6694   & 0.2128 \\
        $\Lambda = 7$ & 0.4168  &  0.4260  &  0.7288   & 0.3387\\
        $\Lambda = 9$ & 0.4232  &  0.4264  &  0.7408   & 0.2788\\ine
        $\mu = 1$ &  0.3884  &  0.4245  &  0.6694  &  0.2056 \\
        $\mu = 3$ &  0.8889  &  0.4306  &  0.6327  &  0.6590 \\
        $\mu = 5$ & 0.8819  &  0.4353  &  0.6387  &  0.6074 \\
        $\mu = 7$ & 0.8782  &  0.4425  &  0.6432  &  0.6345 \\
        $\mu = 9$ & 0.8757  &  0.4443  &  0.6466  &  0.7101 \\
        $\mu = 10$ & 0.8748  &  0.4847  &  0.6481  &  0.6690 \\ine
    \end{tabular}
    \caption{The average absolute relative errors computed from the Equation (\ref{relerr}) quantifying the control errors in correspondence with Figures \ref{fig1} and \ref{fig3}.}
    \label{relative_error}
\end{table}


\section{Conclusion\label{sec_conclusion}}

The nonlinear dynamics control modeling problems of the van der Pol system are tackled by comparing deep learning with traditional deterministic algorithms in this paper. The key idea of this work is to elaborate on the main differences by conducting a comprehensive comparison and benchmark for the recently proposed physics-informed neural networks control with other deterministic algorithms. We first design a benchmark problem for testing the system response for different methods. The desired trajectory and systematic nonlinearity are then changed to check the systematic responses of different controls. The computation burdens are also considered for different methods.

For benchmark analysis, results indicated that all the control theory algorithms exhibit a strong fluctuation which can be interpreted as enlarging the nonlinear inherent van der Pol dynamics with $\mathcal{FF}$ successfully implementing the controls, but the rest fails. The "nonlinearity enlargement" effect is observed to be more obvious for higher-order terms. The PIDOC exhibits stochastic nature, which can be attributed to the nature of deep learning inference, same as reported by Zhai \& Sands \cite{pidoc}. When changing the trajectory amplitudes, an interesting phenomenon is that $\mathcal{FF}$ failed for trajectory convergence when $\Lambda = 1$. Also, a higher control signal implementation density is observed along the horizontal edges of the first order phase portraits, unveiling control theory imposition to van der Pol systems executes stronger controls along the "signal waves' crest and trough." An evidently higher computation burden is observed for PIDOC in comparison to control theory methods. We explain such by the nature of NN learning: the recursive randomization of the NN weights and biases took a longer time than the direct execution of the control signal. For the van der Pol systems with different nonlinearities, it is observed that $\mathcal{FF}$ fails the control when $\mu =10$, whereas PIDOC also failed to implement controls when $\mu \neq 1$, as the controlled schemes were "trapped" into smaller trajectories. The "nonlinearity enlargement effect" for higher-order phase portraits for control theory algorithms. An interesting phenomenon of a vortex-liked structure of the controlled schemes, as originally reported by Zhai \& Sands~\cite{pidoc}, has also been reported for the $\mathcal{C}$ controls. The evidently higher computation time is also reported for PIDOC, the same as what has been reported for different trajectories. For PIDOC, the computation burden is generally reduced with systems of higher nonlinearities. The proposed comparison can guide the future implementation of deep learning-based controller designs and industrial selections.

\section*{\textsc{Data Availability}}

All the data and code will be made publicly available upon acceptance of the manuscript through \url{https://github.com/hanfengzhai/PIDOC}. The {\sf Simulink} file for the deterministic control methods is available upon reasonable requests to the corresponding author. The {\sf Simulink} file was originally published by Cooper \cite{cooper}.



\end{document}